\documentclass[11pt,a4paper]{amsart}
\date{20 June 2011}

\usepackage{latexsym,amsmath,amsfonts,amscd,amssymb, mathrsfs, slashed, amsthm}
\usepackage[english]{babel}

\theoremstyle{plain}  
\newtheorem{theorem}{Theorem}[section]

\newtheorem*{theorem*}{Theorem}

\newtheorem{corollary}[theorem]{Corollary}
\newtheorem{lemma}[theorem]{Lemma}
\newtheorem{proposition}[theorem]{Proposition}

\newtheorem{definition}[theorem]{Definition}
\theoremstyle{remark}
\newtheorem{example}[theorem]{Example}
\newtheorem{remark}[theorem]{Remark}

\newtheorem*{claim*}{Claim}

\numberwithin{equation}{section}

\newcommand{\n}{\noindent}

\renewcommand{\geq}{\geqslant}
\renewcommand{\ge}{\geqslant}

\newcommand{\lra}{\longrightarrow}

\renewcommand{\S}{\slashed{\mathcal{S}}}
\def\lms{\longmapsto}
\def\proof{\noindent\textit{Proof ---}\hspace*{0.2cm}}
\def\qed{\vspace*{-0.1cm} \hfill{$\square$}}
\def\disfrac#1#2{\displaystyle{\frac{#1}{#2}}}

\def\Tr{\mathrm{Tr}}

\def\Ric{\mathrm{Ric}}
\def\Id{\mathrm{Id}}

\begin{document}
\title{Einstein four-manifolds with skew torsion}

\author[A. C. Ferreira]{Ana Cristina Ferreira}
\address{Centro de Matem\'{a}tica \\
Universidade do Minho \\
Campus de Gualtar \\
4710-057 Braga \\
Portugal} \email{anaferreira@math.uminho.pt}


\maketitle

\begin{abstract}
We develop a notion of Einstein manifolds with skew torsion on
compact, orientable Riemannian manifolds of dimension four. We prove
an analogue of the Hitchin-Thorpe inequality and study the case of
equality. We use the link with self-duality to study the moduli
space of 1-instantons on $S^4$ for a family of metrics defined by
Bonneau.
\end{abstract}


\section{Introduction}
\bigskip

Torsion, and in particular skew torsion, has been a topic of interest to both mathematicians and physicists in recent decades. The first attempts to introduce torsion in general relativity go back to the 1920's with the work of \'{E}. Cartan, \cite{Cartan}. More recently, torsion makes its appearance in string theory, where the basic model of type II string theory consists of a Riemannian manifold, a connection with skew torsion, a spinorial field and a dilaton function, \cite{Agricola}.

From the mathematical point of view, skew torsion has played a significant role in the work of Bismut and his local index theorem for non-K\"{a}hler manifolds, \cite{Bismut}. Bismut showed that for any Hermitian manifold, there is a unique connection with skew torsion which preserves both the metric and the complex structure. Nowadays these connections are known as Bismut connections and a Hermitian manifold equipped with such a connection is often referred to as a KT (K\"{a}hler with torsion) manifold. Also, skew torsion is an important feature in Hitchin and Gualtieri's  generalized geometry, \cite{HitchinGen, Gualtieri} where there are natural connections with skew torsion, the exterior derivative of the B-field. In particular, two Bismut connections appear in the characterization of a generalized K\"{a}hler structure.

In this article, we propose a notion of Einstein manifold with skew torsion for a four-manifold. Four dimensions is of particular interest because of the phenomenon of self-duality. We define our notion of Einstein by making use of the decomposition of the curvature operator in terms of the action of $SO(4)$ and making the analogy with the standard Riemannian situation. Motivated by the earlier work of Hitchin and Thorpe, \cite{Hitchin, Thorpe}, we show that an Einstein manifold with skew torsion satisfies a topological constraint, an inequality involving the Euler characteristic and the signature of the manifold: $2 \chi \geq 3 |\tau| $. Manifolds of type $S^1\times S^3$ are well known to satisfy the inequality but they do not carry an Einstein metric. In fact they have a natural structure of flat manifold with skew torsion and we also prove that these are the only manifolds that satisfy the equality $2 \chi = 3 |\tau| $.

Our definition of Einstein with skew torsion depends on the choice of orientation but, as we show, this choice is irrelevant in the compact world. This is seen by establishing a one-to-one correspondence between Einstein manifolds with skew torsion and Einstein-Weyl manifolds and making use of the Gauduchon gauge where the torsion is closed.

An interesting observation is that a connection which is Einstein with skew torsion induces a self-dual connection on the bundle of self-dual forms and, if the manifold is spin, on the bundle of positive half-spinors. In view of this, we use the link with Einstein-Weyl, to present an example of a one-parameter family of U(2)-invariant Einstein metrics with closed skew torsion on the 4-sphere, the triple $(S^4, ds^2, H)$ as defined by G. Bonneau, \cite{Bonneau}. We investigate the moduli space of charge 1 instantons and prove that it is diffeomorphic to that of a metric of constant sectional curvature. The connections with torsion $\pm H$ define two different instantons which we point out to generically define the line of gauge equivalence classes of U(2)-invariant instantons in the moduli space.

\bigskip\bigskip

\section{Metric connections with skew torsion}

Let $(M,g)$ be a Riemannian manifold. Suppose that $\nabla$ is a
connection on the tangent bundle of $M$ and let $T$ be its (1,2)
torsion tensor. If we contract $T$ with the metric we get a (0,3)
tensor which we will still call the torsion of $\nabla$. If $T$ is a
three-form then we say that $\nabla$ is a connection with
skew-symmetric torsion. Given any three-form $H$ on $M$ then there
exists a unique metric connection with skew torsion $H$ defined
explicitly by
$$g(\nabla_X Y, Z) = g(\nabla^g_X Y,Z) +\frac{1}{2} H(X,Y,Z)$$
where $\nabla^g$ is the Levi-Civita connection.

Consider now the triple $(M, g, H)$ and let $\nabla$ be the
connection with skew torsion $H$. If $R^g$ is the Riemannian
curvature tensor and $R^\nabla$ is the curvature tensor associated
with $\nabla$, we can express $R^\nabla$ in terms of $R^g$ and $H$
as follows: for every four vector fields $X, Y, Z, W$, we have

\begin{eqnarray}\label{Eq: curvature}
 R^\nabla(X,Y,Z,W) & = & R^g(X,Y,Z,W)  \\
& & + \disfrac{1}{4} g(H(X,W),H(Y,Z)) - \disfrac{1}{4} g(H(Y,W), H(X,Z)) \vspace*{2mm} \nonumber \\
& & - \disfrac{1}{2}(\nabla^g_X H)(Y,Z,W) + \disfrac{1}{2}
(\nabla^g_Y H)(X,Z,W). \nonumber
\end{eqnarray}

The tensor $R^\nabla$ will not have the same symmetries as $R^g$, as
can be expected. For example, the analogue of the Bianchi identity
is given by

\begin{equation}\label{Eq: Bianchi}
 R^\nabla(X,Y,Z,W)+R^\nabla(Y,Z,X,W) +
R^\nabla(Z,X,Y,W) =     \end{equation} $$\hspace*{1cm} - dH(X,Y,Z,W)
- (\nabla_W^g H)(X,Y,Z)+\disfrac{1}{2} \mathop{\sigma}_{XYZ}
g(H(X,Y),H(Z,W)),$$

where $\displaystyle{\mathop{\sigma}_{XYZ}}$ denotes the cyclic sum
over $X,Y,Z$.

\medskip

Also, if $\nabla^-$ is the metric connection with torsion $-H$, we
have

\begin{equation}\label{Eq: swap-general}
R^\nabla(X,Y,Z,W) = R^{\nabla^-}(Z,W,X,Y) - \disfrac{1}{2}~ d
H(X,Y,Z,W)
\end{equation}

In particular if $H$ is closed, we obtain

\begin{equation}\label{Eq: swap-particular}
R^\nabla(X,Y,Z,W) = R^{\nabla^-}(Z,W,X,Y).
\end{equation}

Suppose now that $M$ is orientable. If $n = \dim M$, let
$\{e_i\}_{i=1}^n$ denote a positively oriented orthonormal frame of
$TM$ and $\{e^i\}_{i=1}^n$ its dual frame. Also, let $\Ric^g$ be the
usual Ricci tensor and $\Ric^\nabla$ be the Ricci tensor with
respect to $\nabla$. Given any pair of vector fields $X,Y$ we have

\begin{equation}\label{Eq: Ricci}
\Ric^\nabla(X,Y)=\Ric^g(X,Y) - \disfrac{1}{4} \sum_i
g(T(X,e_i),T(Y,e_i))-\disfrac{1}{2} d^* H(X,Y).
\end{equation}

Notice that, unlike $\Ric^g$, $\Ric^\nabla$ has a non-vanishing
anti-symmetric part.

For the scalar curvature, the relation is

\begin{equation}\label{Eq: scalar}
s^\nabla-s^g=-\disfrac{3}{2}\Vert H \Vert^2
\end{equation}

These or similar identities are already available in the literature.
The original proof of \ref{Eq: swap-general}-\ref{Eq:
swap-particular} was done in \cite{Bismut}. For the remaining ones,
see \cite{IvanPap}.

\bigskip

\section{Decomposition of the Riemann tensor}

\bigskip

We now restrict our attention to manifolds of dimension four. We
will also be assuming compactness and orientability. Recall that the
bundle of two-forms splits as $\Lambda^2 = \Lambda_+\oplus \Lambda_-
$, where $\Lambda_+$ and $\Lambda_-$ are the bundles of self-dual
and antiself-dual forms, respectively.

Consider $(M,g, H)$ and notice that in four dimensions, the star
operator also allows us to see the torsion $H$ as a one-form. Let
$h=*H$ and call $h$ the torsion 1-form. One of the features of $h$
is that it provides us with a simple way of writing the expression
for the Ricci tensor.

\medskip

\begin{proposition}\label{Prop: Ricci-4dim}
On a four-dimensional manifold, the Ricci tensor for a connection
$\nabla$ with skew torsion $H$ can be written as
$$\Ric^\nabla = \Ric^g- \disfrac{1}{2} \Vert h \Vert^2 g+ \disfrac{1}{2} h \otimes h -
\disfrac{1}{2} *dh,$$ where $g$ is the metric tensor.
\end{proposition}

\proof Direct computation using the orthonormal frame $\{e_i\}$.

\qed

\bigskip

The curvature tensor $\,R^\nabla\,$ of a connection with skew
torsion lives in $\Lambda^2 \otimes\Lambda^2$. Using the metric, we
can see $R^\nabla$ as a map $\mathcal{R}^\nabla:
\Lambda^2\lra\Lambda^2$, called the curvature operator, which is
given by the prescription
$$g\left(\mathcal{R}^\nabla\left(X\wedge Y\right), Z\wedge W\right) = R^\nabla\left(X,Y,Z,W\right). $$

We are going to work out the decomposition of $\mathcal{R}^\nabla$
in terms of the splitting $\Lambda^2 = \Lambda_+ \oplus \Lambda_-.$
First let us recall briefly what happens in the usual Riemannian
situation. The symmetries of $R^g$ mean that this is an element of
$S^2\Lambda^2$, which can be decomposed as follows, \cite{Besse}.

$$\displaystyle{\mathcal{R}^g = \left(
\begin{array}{c|c}
& \\ W^+ + \disfrac{s}{12} \Id  & Z \\  & \\  \hline & \\ Z^t &
W^-+\disfrac{s}{12}\Id\\ &
\end{array}
\right)}$$ where $s$ is the scalar curvature, $W^+$ and $W^-$ are
the self dual and antiself-dual parts of the Weyl tensor, and
$Z$ is the trace-free part of the Ricci tensor, i.e.,
$Z=\Ric^g-\frac{s}{4} g$.

\medskip

For a connection with skew torsion the situation is slightly more
complicated, since $R^\nabla$ has a non-vanishing part in
$\Lambda^2(\Lambda^2)$.
\medskip

\begin{theorem}\label{Theo: riemanndecomposition}
For a metric connection with skew torsion $\nabla$, we can decompose
the Riemann curvature map $\mathcal{R}^\nabla$ in terms of self-dual
and antiself-dual blocks as:

$$\displaystyle{\small{\mathcal{R}^\nabla = \left(
\begin{array}{c|c}
& \\ W^+ + \left(\disfrac{s^\nabla}{12} - \disfrac{*dH}{4}
\right)\Id + \disfrac{1}{4}\left(d^* H\right)_+  & Z^\nabla +
S \left( \nabla\! *\! H \right)+ \disfrac{*dH}{4}  g  \\  & \\  \hline & \\
 \left( Z^\nabla - S \left(\nabla\! *\! H \right)+ \disfrac{*dH}{4}  g  \right) ^\dag & W^-+\left( \disfrac{s^\nabla}{12}+ \disfrac{*dH}{4} \right) \Id - \disfrac{1}{4}\left(d^* H\right)_-\\
&
\end{array}
\right)}}$$

\noindent where $S$ denotes the symmetrization of a tensor and
$\dag$ the transpose of a matrix, $Z^\nabla$ is the symmetric
trace-free part of $\Ric^\nabla$, and $\left(d^* H\right)_+$ and
$\left(d^* H\right)_-$ are the self-dual and antiself-dual part of
$d^* H$, respectively.
\end{theorem}

\proof We start with the upper left block and call it A. The best
way to see what the entries are is to do an example. We will be
using the convention $R_{ijkl}$ for $R\left(e_i,e_j,e_k,e_l\right)$.
Write $R = R^g + \overline{R}$ and recall equation \ref{Eq:
curvature}. Take the first diagonal entry, this is given by
$$A_{11} = \disfrac{1}{2}\left(R_{1212}+R_{1234}+R_{3412}+R_{3434}\right).$$
We only need to worry about the $\overline{R}$ component. We can
easily see that
$$\overline{R}_{1212}+\overline{R}_{3424} = -\frac{1}{4}\left(({H_{12}^r})^2 +  ({H_{34}^r})^2\right) = -\disfrac{1}{4} \Vert H \Vert^2$$
and that $$\overline{R}_{1234}+\overline{R}_{3412} =
\frac{1}{2}\left(H_{14}^r H_{23}^r  -H_{13}^r H_{24}^r \right)-
\frac{1}{2}\left( dH \right)_{1234}$$ and since, $H_{ij}^rH_{kl}^r$
vanishes if $i,j,k,l$ are all distinct, we get
$$\overline{A}_{11} = -\frac{\Vert H \Vert^2}{8}  - \frac{*dH}{4}$$ and the same holds for the other diagonal entries.
Consider the off-diagonal entries now. Taking $A_{12}$, for
instance, we get that
$$\overline{A}_{12} = - \frac{1}{8}\left( H_{12}^r H_{13}^r - H_{12}^r H_{24}^r + H_{34}^r H_{13}^r - H_{34}^rH_{24}^r \right) + \frac{1}{4}
\left( \left(d^* H\right)_{14} + \left(d^* H\right)_{23} \right).
$$  The quadratic part of the expression vanishes, so we
obtain
$$\overline{A}_{12} = \frac{1}{4}
\left( \left(d^* H\right)_{14} + \left(d^* H\right)_{23} \right)
$$ and the other off-diagonal entries are analogous. Then, clearly,
we have
$$A  =  W^+ + \left( \frac{s}{12} - \frac{\Vert H \Vert^2}{8} -
\frac{*dH}{4} \right) \Id + \frac{\left(d^* H\right)_+}{4}
 $$ and since $s^\nabla = s^g - \frac{3}{2}\Vert H \Vert^2$, by
equation \ref{Eq: scalar}, we get the desired expression for $A$. If
$D$ is the lower right block then the arguments are perfectly
similar to the ones for $A$. Consider now the upper right block,
$B$. Let us start with the diagonal entries. Take
$$\overline{B}_{11}=\frac{1}{2}\left(\overline{R}_{1212}+\overline{R}_{1234}-\overline{R}_{3412}-\overline{R}_{3434}\right).$$
We have
$$\begin{array}{lcl}
\overline{R}_{1212}-\overline{R}_{3434} & = & -\frac{1}{4}
\left((H_{123})^2+ (H_{124})^2-(H_{134})^2-(H_{234})^2\right) \vspace*{2mm}\\
 \overline{R}_{1234}-\overline{R}_{3412} & = &
-\frac{1}{2}\left(\left(\nabla^g_{1}
H\right)_{234}-\left(\nabla^g_{2} H\right)_{134}- \left(\nabla^g_{3}
H\right)_{412}+ \left(\nabla^g_{4} H\right)_{312}\right)
\end{array}$$
and we will now write this in terms of $h=*H$, since it makes the
calculations easier. We then have that
$$\overline{B}_{11}= \frac{1}{8} \left((h_{1})^2+(h_{2})^2-(h_{3})^2-(h_{4})^2\right) +
\frac{1}{2}\left(\left(\nabla^g_{1} h\right)_{1}+\left(\nabla^g_{2}
h\right)_{2}- \left(\nabla^g_{3} h\right)_{3} - \left(\nabla^g_{4}
h\right)_{4}\right)$$ and analogously for $\overline{B}_{22}$ and
$\overline{B}_{33}$. Consider now
$$\overline{B}_{12} = \frac{1}{2}\left( \overline{R}_{1312}+ \overline{R}_{1334}+\overline{R}_{2412}+\overline{R}_{2434}
\right)$$ and we see that
$$\overline{B}_{12} = \frac{1}{4}\left( H_{123}H_{234} - H_{124}H_{134} + \left(\nabla^g_1 H\right)_{132}
-\left(\nabla^g_3 H\right)_{314}-\left(\nabla^g_2 H\right)_{241}+
\left(\nabla^g_4 H\right)_{423} \right)$$ and, rewriting in terms of
$h$, we get $$\overline{B}_{12}=\frac{1}{4}\left(h_2h_3-h_1h_4 +
\left(\nabla^g_2 h\right)_3 +\left(\nabla^g_3 h\right)_2
-\left(\nabla^g_1 h\right)_4 - \left(\nabla^g_4 h\right)_1 \right)$$
and we have similar results for the other entries. We wish to
express $B$ in terms of symmetric trace-free 2-tensors, so we need
to choose the right isomorphism between $\Lambda_+\otimes\Lambda_-$
and $S_0^2$. This can be found in \cite{Besse}, is called the Ricci
contraction and works as follows: if we consider 2-forms as matrices
then $-\varphi$ is given by standard matrix multiplication. For
example, the form $e^1\wedge e^2$ corresponds to the $4\times 4$
matrix $M$ such that $M_{21} = 1, M_{12} = -1$ and $M_{ij}=0$
elsewhere.

We can now conclude that $$\overline{B}= \frac{1}{2} \left(h\otimes
h - \frac{1}{4} \Vert h \Vert^2 g \right) + S\left(\nabla^g h \right)+
\frac{d^*h}{4}g$$ where $S$ denotes the symmetrization of the
tensor. Observe the following two lemmas, which can be proved by
simple local calculations.

\begin{lemma}\label{Lem: same-cov-der}
If $\nabla$ is the metric connection with skew torsion $H$, and $h=
* H$, then $\nabla \! h = \nabla^g\! h$.
\end{lemma}

\begin{lemma}
The trace-free symmetric part of the Ricci tensor $\Ric^\nabla$,
denoted by $Z^\nabla$, is given by:
$$Z^\nabla = Z^g + \disfrac{1}{2} h\otimes h - \disfrac{1}{8} \Vert h \Vert^2 g.$$
\end{lemma}

\medskip

\n Finally, we get that $$B=Z^\nabla+S\left(\nabla\! *\! H \right)+
\frac{*dH}{4}g.$$ If $C$ is the remaining block, by means of
equation \ref{Eq: swap-general} and noticing that we always have two
repeated indices,  $C$ is the transpose of $B$ when replacing $H$ by
$-H$.

\qed

\bigskip\bigskip

\section{Einstein metrics with skew torsion}\label{Sec: Einstein-skew}

\bigskip

The above decomposition of the Riemann tensor of a connection with
skew torsion is our main motivation for the following definition,
recalling also that in standard Riemannian geometry, a manifold
$(M,g)$ is said to be Einstein if $Z^g=0$.

\medskip

\begin{definition}\label{def: Einstein-metric}
Given an oriented Riemannian four-manifold $(M, g, H)$, we say that
$g$ is an Einstein metric with skew torsion, if $$Z^\nabla + S\left(
\nabla
*H \right) + *\frac{dH}{4}g  = 0$$ where $\nabla$ is the metric
connection with skew torsion $H$.
\end{definition}

We remark that the standard notion of Einstein metric is equivalent
to having the induced Levi-Civita connections on $\Lambda_+$ and
$\Lambda_-$ self-dual and anti-self-dual, respectively. Our
definition of Einstein with skew torsion simply adapts this, but we
usually do not have both statements in our situation. Here we have
chosen that $\nabla$ on $\Lambda_+$ be self-dual. We see will later,
in corollary \ref{Cor: XXX}, that for a compact manifold this choice
does not constitute a problem.

\medskip

\begin{example}
The very basic example is the one of the Lie group $S^1\times S^3$,
with one of the two flat connections given by left or right
trivialization of the tangent bundle.
\end{example}

\medskip

\begin{example}

Recall that the equations of type II string theory may be
geometrically described as a tuple $(M,g,H,\phi,\psi)$ consisting of
a manifold $M$ with a Riemannian metric $g$, a three-form $H$, a
so-called dilaton function $\phi$ and a spinor field $\psi$
satisfying the following system of equations, \cite{Agricola},
$$\begin{array}{ll}
\Ric^\nabla +\frac{1}{2} d^*H + 2 \nabla^g d\phi =0 \qquad&
(\nabla^g_X
+\frac{1}{4} X\lrcorner H)\psi = 0 \vspace*{1mm}\\
d^*(e^{-2\phi}H) = 0 & (2d\phi -H)\psi = 0
\end{array}$$
where $\nabla=\nabla^g+\frac{1}{2}H$. Suppose  $2d\phi =
*H$, then the first equation implies
$$S(\Ric^\nabla) + \nabla\! *\! H = 0$$ since $\nabla \! * \! H$
is the Hessian of $2\phi$ and is therefore symmetric. Hence the
trace-free part satisfies definition \ref{def: Einstein-metric}.
\end{example}

\medskip

We have an interesting property in the compact case if the torsion
is closed.

\begin{proposition}\label{Prop: killing}
If $M$ is compact and  $dH=0$, the Einstein equations with skew
torsion imply that the vector field $X$ defined by $i_X \omega_g =
H$, where $\omega_g$ is the volume form, is a Killing field.
\end{proposition}

\proof It suffices to prove that
$$\int_M ||S(\nabla^g h)||^2\omega_g = 0,$$ where $h$ is the one-form dual to $X$. We can write $S(\nabla^g h)$ as $\nabla^g h -\frac{1}{2}d h$ and
using the Einstein condition with skew torsion also as
$-(Z^g-\frac{1}{8}||h||^2g+\frac{1}{2}h \otimes h)$. Observing that
the decomposition $T^*\otimes T^* = S^2 T^* \oplus \Lambda^2 T^*$ is
orthogonal we get that
$$\int_M ||S(\nabla^g h)||^2 \omega_g = \int_M -\left( Z^g-\frac{1}{8}||h||^2g+\frac{1}{2}h\otimes h,\nabla^g h\right)\omega_g$$
where the round brackets denote here the inner product of tensors,
for convenience. It is easy to see that $(g,\nabla^g h)= 0$, since
$d^* h=0$; therefore we are left with
$$\int_M -\left(\Ric^g + \frac{1}{2} h\otimes h, \nabla^g h\right)\omega_g.$$
Recall that the divergence of a two-symmetric tensor is given by
$({\nabla^g})^*$, the formal adjoint of $\nabla^g$. Recall also that
by contracting the differential Bianchi identity, we get that the
divergence of $\Ric^g$ is $-\frac{1}{2}ds^g$. Then we can write the
 integral as $$\int_M \left( \frac{1}{2} ds^g-\frac{1}{2}{\nabla^g}^*
(h \otimes h), h \right)\omega_g.$$ The idea now is to write this as
a divergence; since $d^*h =0$, $(ds^g,h)= -{\nabla^g}^* (s^g\, h)$
and $({\nabla^g}^* (h \otimes h), h)=\frac{1}{2}{\nabla^g}^*(||h||^2
h)$. Finally we obtain that
$$\int_M  {\nabla^g}^*\left(- \frac{1}{2} s h + \frac{1}{4} ||h||^2 h \right) \omega_g = 0.$$

\qed

\medskip
The above result also derives from another interpretation of the
Einstein equations with skew torsion which is related to conformal
invariance and was originally proved in this context, \cite{Tod}.

\medskip

\begin{corollary}\label{Cor: XXX}
On a compact four-manifold, an Einstein metric with closed skew
torsion $H$ satisfies the equation $Z^\nabla = 0$, where $Z^\nabla$
is the trace-free part of the Ricci tensor.
\end{corollary}

\proof Recall that a vector field $X$ is said to be Killing if the
symmetric part of $\nabla^g X$ vanishes. If $X$ is the metric dual
of $h$, where $h= *H$, then
$$\nabla^g X = \nabla^g h$$ and by means of Lemma \ref{Lem: same-cov-der} $$\nabla^g h= \nabla h.$$
Using proposition \ref{Prop: killing}, then definition \ref{def:
Einstein-metric} gives $Z^\nabla = 0$.

\qed
\medskip

\begin{remark}\label{Rem: plus-minus-H}
It is clear, using the corollary above and looking at the expression
of
$$Z^\nabla = Z^g + \frac{1}{2} * H \otimes * H - \frac{1}{8}||H||^2\,
g,$$ that if $(M, g, H)$ is a compact Einstein manifold with closed
skew torsion then so is $(M, g, -H)$.
\end{remark}

\bigskip\bigskip

\subsection{An inequality}\label{Subsec: inequality}
\bigskip

Our definition of Einstein metric with skew torsion implies that
$\Lambda_+$ has a self-dual connection. This means that
$\mathrm{Tr}(R\wedge R) = f \omega_g$, where $f$ is a non-negative
function, and hence the first Pontryagin class of $\Lambda_+$ is
non-negative. This implies a topological constraint on a compact
four-manifold that generalizes the Hitchin-Thorpe inequality,
\cite{Besse, Hitchin}, which states that if $M$ is a compact
oriented Einstein manifold of dimension 4, then the Euler
characteristic $\chi(M)$ and the signature $\tau(M)$ satisfy the
inequality $\chi(M)\geq \frac{3}{2}\vert \tau(M) \vert.$

\medskip

We have a similar result for connections with skew torsion.

\medskip
\begin{theorem}\label{Theo: inequality}
Let $(M,g, H)$ be a compact, oriented, four-dimensional Riemannian
manifold, equipped with a metric connection with skew-symmetric
torsion $H$, such that $Z^\nabla + S\left(\nabla\! *\! H +
\frac{*dH}{4}g \right) = 0$, then
$$\chi(M)\geq \disfrac{3}{2}\left|\tau(M)\right|.$$
\end{theorem}

\proof We use the formulas discussed and proved in \cite{Besse2}.
Both the Euler characteristic and the signature can be written in
terms of the curvature operator $\mathcal{R}$ as
$$\chi(M)= \frac{1}{8\pi^2}\int \Tr(\mathcal{*R*R})\omega_g$$
$$\tau(M) = \frac{1}{12\pi^2}\int \Tr(\mathcal{R*R})\omega_g$$
where $*$ is the Hodge star operator and $\omega_g$ the volume form
with respect to the metric and the chosen orientation. Recall from
the proof of theorem \ref{Theo: riemanndecomposition} that
$\mathcal{R}$ is given in blocks by
$$\mathcal{R}=\left(
\begin{array}{cc}
A & B\\
C & D
\end{array}
 \right).$$
Clearly, we have that $*A=A$, $*D=-D$ and since $B=0$ we get
$\Tr(\mathcal{*R*R})=\Tr(A^2+D^2)$ and
$\Tr(\mathcal{R*R})=\Tr(A^2-D^2)$. Observe now that
$$\begin{array}{c}\Tr(\mathcal{*R*R})= \Tr(A^2+D^2)\ge \Tr(A^2-D^2) = \Tr(\mathcal{R*R})\\
\Tr(\mathcal{*R*R})= \Tr(A^2+D^2)\ge \Tr(D^2-A^2) =
-\Tr(\mathcal{R*R})\end{array}$$ which gives two inequalities
$\chi(M)\ge \frac{3}{2}\tau(M)$ and $\chi(M)\ge
-\frac{3}{2}\tau(M)$, and combining these two we get the desired
inequality.
 \qed

\bigskip

Our next goal is to determine for which Einstein metrics with skew
torsion the equality is attained. First we need to look at the
conformal class of such a metric.

\bigskip\bigskip

\subsection{Conformal invariance}
\bigskip

We now introduce the notions of Weyl structure and Einstein-Weyl
manifold, \cite{CalPed}.

\begin{definition}
Let $M$ be a manifold with conformal structure $[g]$, i.e., an
equivalence class of metrics such that $\tilde{g} \simeq g$ if
$\tilde{g} = e^f g$, where $f:M\lra \mathbb{R}$ is a smooth
function. A Weyl connection is a torsion-free affine connection $D$
such that for any representative of the metric $g$ there exists a
one-form $\omega$ such that $Dg=\omega\otimes g$. A Weyl manifold is
a manifold equipped with a conformal structure and a compatible Weyl
connection. The Weyl structure is said to be closed (resp. exact) if
(any) $\omega$ is closed (resp. exact).
\end{definition}

We note that the notions of closed and exact Weyl structures are
well defined. If $\omega$ is the one-form associated to $g$ and
$\tilde{\omega}$ is the one-form associated to $\tilde{g}= e^f g$,
then $\tilde{\omega} = \omega + d f$.

\medskip

\begin{definition}
A Weyl manifold is said to be Einstein-Weyl if the trace-free
symmetric part of the Ricci tensor $S_0(\Ric^D)$ vanishes.
\end{definition}

\medskip

The following formulas, \cite{PedSwa2}, are simple but extremely
useful calculations:

The Weyl connection $D$ with one-form $\omega$ is given explicitly
by
\begin{equation}
D=\nabla^g_X Y -\frac{1}{2} \omega(X)Y - \frac{1}{2}\omega(Y)X +
\frac{1}{2}g(X,Y)\omega^\sharp \label{eq: Weyl-connection}
\end{equation} where $\omega^\sharp$ denotes the vector field dual
to $\omega$. The symmetric part of its Ricci tensor is equal to
\begin{equation}
S(\Ric^D) = \Ric^g - \frac{1}{2} (\Vert \omega \Vert^2 g -
\omega\otimes\omega ) + S(\nabla^g \omega) - \frac{1}{2}(d^*\omega)
g.\label{eq: Ricci-Weyl} \end{equation}

\medskip

This immediately yields,

\begin{theorem}\label{theo: Einstein-Weyl}
Let $(M,g,H)$ be a four-dimensional Einstein manifold with skew
torsion. Then if $\omega = *H$, the torsion-free connection $D$ such
that $D\omega = \omega\otimes g$ is an Einstein-Weyl connection.
Conversely, given an Einstein-Weyl manifold, each metric in the
conformal class defines, with $H= -*\omega$, an Einstein manifold
with skew torsion.
\end{theorem}

\proof Suppose $(M,[g])$ is Einstein-Weyl. Take a representative of
the metric $g$ and its associated one-form $\omega$. The connection
defined by equation $\ref{eq: Weyl-connection}$ has scalar curvature
$$s^D=s^g - \frac{3}{2}\Vert \omega \Vert^2-3 d^*\omega.$$ Therefore, using also equation \ref{eq: Ricci-Weyl}, the trace-free
symmetric Ricci tensor is equal to
$$S_0(\Ric^D)=\Ric^g + \frac{1}{2} \omega\otimes \omega - \frac{1}{8} \Vert \omega \Vert^2 g + S(\nabla^g \omega) +
\frac{1}{4} (d^* \omega) g.$$ Now take the metric connection with
skew torsion $H=-*\omega$. Then clearly $(M,g,H)$ is Einstein with
skew torsion. The converse is perfectly analogous.

\qed

\medskip

As an immediate corollary of this, we get that the Einstein
equations with skew torsion are conformally invariant, that is, if
the metric $g$ is Einstein with skew torsion, then so are all
metrics in the conformal class of $g$, if we transform the torsion
appropriately.

Notice again that, unlike in string theory and Einstein-Weyl
geometry, definition \ref{def: Einstein-metric} does not work in any
dimension except four. Indeed, it is crucial that $*H$ is a
one-form.

Still in the context of conformal invariance we have the following
important fact:  given a metric $g$ on a compact manifold and a
one-form $\omega$, there is a unique (up to a constant) metric
$\tilde{g} = e^f g$ for some smooth function $f$, such that the
one-form $\tilde{\omega} = \omega + d f$ is co-closed with respect
to $\tilde{g}$. This metric is of particular importance in Hermitian
geometry and it is known in the literature as the Gauduchon gauge,
\cite{Gauduchon1}. We, then, have the following,

\medskip

\begin{corollary}\label{Cor: YYY}
If $(M, g, H)$ is a compact Einstein manifold with skew torsion then
there exists a function $f$ on $M$ such that $(M, e^f g, e^f (H - *df))$ is Einstein with closed skew torsion.
\end{corollary}

\medskip

The above corollary together with corollary \ref{Cor: XXX} implies
that our definition of Einstein metrics with skew torsion is
independent of orientation in the case of compact manifolds.

It should also be mentioned here that a generalization of the
Hitchin-Thorpe inequality for Einstein-Weyl manifolds was proved in
\cite{PedPooSwa}.

\bigskip\bigskip

\subsection{The equality}\label{Subsec: equality}
\bigskip

As mentioned in subsection \ref{Subsec: inequality} we are
interested in the case where equality is achieved. The usual
Riemannian situation was studied by N. Hitchin \cite{Besse,
Hitchin}, who proved that if $M$ is a compact oriented
four-dimensional Einstein manifold and the Euler characteristic
$\chi(M)$ and the signature $\tau(M)$ satisfy
$$\chi(M)=\frac{3}{2}\vert \tau(M) \vert $$ then the Ricci curvature
vanishes, and $M$ is either flat or its universal cover is a $K3$
surface. In that case, $M$ is either a $K3$ surface itself
$(\pi_1(M)=1)$, or an Enriques surface $(\pi_1(M)=\mathbb{Z}_2)$, or
the quotient of an Enriques surface by a free antiholomorphic
involution $(\pi_1(M)=\mathbb{Z}_2\times \mathbb{Z}_2)$ with the
metric induced from a Calabi-Yau metric on $K3$.

\medskip

In the following we investigate what happens when equality holds in
our setting of connections with skew torsion. Given the link with
Einstein-Weyl geometry, it is not surprising that a classification
has been achieved for the four-dimensional case with closed Weyl
structure, \cite{Gauduchon2}. This is somewhat similar to what we want,
but the arguments rely on  twistor theory which we
want to avoid, \cite{Gauduchon2}. Instead we will keep to the
language of Riemannian geometry.

\medskip

\begin{theorem}\label{Theo: einstein-skew}
Let $(M,g,H)$ be a Riemannian compact, oriented four-manifold M
which is an Einstein manifold with skew torsion satisfying the
equality
$$\chi(M)= \frac{3}{2} |\tau(M)|.$$ The either $M$ is Einstein or its
universal cover is isometric to $\mathbb{R}\times S^3$.
\end{theorem}

\medskip

\begin{remark}
As mentioned before $M=S^1\times S^3$ is a compact solution of the
Einstein equations with skew torsion. Also observe that $S^1\times
S^3$ is not an Einstein manifold in the usual sense. Since
$\chi(M)=0$, if $(M,g)$ was Einstein then we would have
$$\chi(M)=\disfrac{1}{8\pi^2}\int_M \left( \frac{s^2}{24}+||W||^2\right) \omega_g$$ which would mean that both the scalar curvature and the Weyl tensor vanish. Therefore $S^1\times S^3$ would be flat with respect to the Levi-Civita connection which is a contradiction.
\end{remark}

\medskip

\n\textit{Proof of theorem \ref{Theo: einstein-skew} ---} From
corollary \ref{Cor: YYY}, we can assume that $(M,g,H)$ is such that
$dH =0$. Let $h$ be the torsion one-form. Suppose, without loss of
generality, that $\chi(M)=-\frac{3}{2}\tau(M)$. Then $\Tr(A^2)=0$
(recall the proof of theorem \ref{Theo: inequality}) and from the
decomposition of $\Lambda^2\otimes \Lambda^2$ into irreducible
$SO(4)$-components, we get
$$\Vert W^+\Vert^2 = \Vert s^\nabla\Vert^2 = \Vert (d^*H)_+\Vert^2 = 0.$$
Then, in particular, $*dh$ is anti-self-dual, and we have
$$-\Vert dh\Vert^2 = \int_M *dh \wedge *dh = \int_M dh\wedge dh = \int_M d(h\wedge dh)$$
and so $dh$ vanishes by Stokes theorem.  Recall, from lemma
\ref{Prop: killing}, that if $X$ is the dual of $h$ via the metric
$g$, then $X$ is a Killing field. Combining these two facts we
conclude that $\nabla^g X = 0$. Then either $X=0$ and we are in the
Einstein situation or otherwise $X$ is a nowhere vanishing parallel vector field.
In this case we have a reduction of the holonomy group and, by means of the de Rham decomposition theorem, $M$ splits locally as a Riemannian product $\mathbb{R}\times N$.
Since $\Ric^\nabla = 0$ then $$\Ric^g = \frac{1}{2}\Vert
h\Vert^2 g  -\frac{1}{2} h\otimes h.$$
Observing that $TN$ is
the orthogonal complement of $\{X\}$, we conclude that $N$ is Einstein with positive Ricci curvature. Hence, since $N$ is of dimension 3, it is of positive sectional curvature. Therefore $M$ is
locally isometric to $\mathbb{R}\times S^3$, the metric splits as a
product and the three-form is the pull-back of a three form in $N$, using the inclusion.

\qed

\begin{remark}
A natural question to ask is which compact Hermitian four-manifolds equipped with the Bismut are Einstein in the sense of our definition
\ref{def: Einstein-metric}. The answer does not give new examples of such manifolds. We can prove that the Lee form is parallel and then repeat the steps of the proof above,
\cite{Ferr}.

\end{remark}

\medskip

Given what was presented here so far, a natural question to ask is if there are other instances of compact Einstein manifolds besides manifolds of  type $S^1\times S^3$. The answer is yes, and given theorem \ref{theo: Einstein-Weyl}, a good source of examples is that of Einstein-Weyl geometry.

We can find a classification of
four-dimensional Einstein-Weyl manifolds whose symmetry group is at
least four dimensional in \cite{Madsen}. This article has two errors in the case
of $U(2)$-invariant structures which were pointed out by G. Bonneau
in \cite{Bonneau} who also offers a simpler description of the
metrics in the Gauduchon gauge. We can summarize the results for the compact orientable case as
follows:

If $(M, g, H)$ is a compact orientable four-dimensional manifold
which is Einstein with closed skew torsion and whose symmetry group
is at least four-dimensional, then we have one of the following
possibilities:
\begin{itemize}
\item if $*H$ is exact then $M$ is Einstein,
\item if $*H$ is closed but not exact then $M$ is finitely covered
by $S^1\times S^3$ with its standard flat structure,
\item if $*H$ is not closed then the symmetry group is \\
- $S^1\times
SO(3)$ in which case $M$ is $S^4$, $S^1\times S^3$, $S^1\times_{(-1,-1)} S^3$,
$S^2\times S^2$ or $S^2\times_{(-1,-1)} S^2$,\\
-  $U(2)$ in which case $M$ is $S^4$,
$\mathbb{C}P^2$ or $\mathbb{C}P^2\#\overline{\mathbb{C}P}^2$.
\end{itemize}
Also, for each of the listed manifolds there is, in fact,  an Einstein structure with skew torsion.

\bigskip \bigskip

\section{Instantons}

As remarked in section \ref{Sec: Einstein-skew}, if we have an Einstein manifold with skew torsion $(M,g, H)$ and $\nabla$ is the metric connection with torsion $H$ then the induced connection on $\Lambda^+$ is self-dual. In the particular case where $M$ is a spin manifold, the induced connection on $\S^+$, the bundle of positive half-spinors, is also self-dual. Self-dual connections are also called instantons.

By remark \ref{Rem: plus-minus-H}, in the compact case with closed $H$, we will have two different instantons. A question that arises here is whether or not such two instantons are always gauge equivalent. Note that if $\nabla^+$ and $\nabla^-$ represent the two induced connection with torsion $H$ and $-H$ respectively, the Yang-Mills density is the same, namely
$$\left\vert F^{\nabla^\pm} \right\vert = \left( \left\vert W^+\right\vert^2 +  \left\vert  \frac{s^{\nabla^\pm}}{12} \Id\right\vert^2 + \left\vert \frac{(d^*H)_+}{2}  \right\vert^2 \right) \mathrm{dvol}.$$
We will now point out certain features of these concepts for a particular example.

\subsection{An example on $S^4$}

\bigskip

Consider the family of $U(2)$-invariant
metrics on $S^4$, mentioned in subsection \ref{Subsec: equality}.
This is a one parameter family of Einstein metrics with skew
torsion, presented explicitly in \cite{Bonneau} in diagonal form by

\begin{equation}\label{Eq: Bonneau-metrics}
ds^2 = \frac{2}{\Gamma} \left[
\frac{k-x}{\Omega^2(x)(1+x^2)^2}(dx)^2 +
\frac{k-x}{1+x^2}[(\sigma^1)^2+(\sigma^2)^2] +
\frac{\Omega^2(x)}{k-x}(\sigma^3)^2  \right]
\end{equation}

\begin{equation}\label{Eq: Bonneau-forms}
H =  2\disfrac{k-x}{(1+x^2)^2} dx\wedge \sigma^1 \wedge \sigma^2
\end{equation}

 where $x \in (-\infty, k)$ is a
coordinate,  $\{\sigma^i, i=1,2,3\}$ is a basis of left-invariant
forms such that \hbox{$d\sigma^i = \frac{1}{2} \epsilon_{ijk}
\sigma^j \wedge \sigma^k$},
$$\Omega^2(x) = 1+ n (x^2-1-2kx)\left( \frac{\pi}{2}  +
\arctan(x)\right)+ n(x-2k),$$ $\Gamma$ is a positive homothetic
parameter, $k$ is a free parameter and $n$ is such that $$n =
\disfrac{1}{k + (1+k^2) \left(\frac{\pi}{2} + \arctan(k)\right)}.$$
Since $\Gamma$ is simply a homothetic parameter, we can take it to
be $\Gamma = 2$, for simplicity of calculations. For ease of
notation, we will be writing
$$ds^2 = a^2\, dx^2 + b^2\, [(\sigma^1)^2 + (\sigma^2)^2] + c^2\, (\sigma^3)^2
$$ where
$$a^2= \frac{k-x}{\Omega^2(x)(1+x^2)^2}, \quad b^2=
\frac{k-x}{1+x^2},\quad c^2= \frac{\Omega^2(x)}{k-x}.$$

\medskip

We are interested in studying charge 1 instantons for $(S^4, ds^2)$ under the gauge group $SU(2)$.
Recall that such a moduli space for a round metric, i.e. a metric of constant sectional curvature, has been studied intensively and is well understood, \cite{AHS, DonKron}. It is a five dimensional manifold which is diffeomorphic to hyperbolic space $\mathbb{H}^5$. A natural question for us is whether the moduli space for each of the Bonneau metrics is smooth and if so what type of manifold is it? The answer to these queries can be
obtained with the help of the following theorem.

\medskip

\begin{theorem}[Buchdahl, \cite{Buch}]\label{Theo: buch}
Let $X$ be a compact complex surface biholomorphic to a blow-up of
$\mathbb{C}P^2$ $n$ times, and $L_{\infty}\subset X$ be a rational
curve with self-intersection $+1$. Let $Y$ be a smooth four-manifold
diffeomorphic to $n\mathbb{C}P^2$ obtained by collapsing
$L_{\infty}$ to a point $y_{\infty}$ and reversing the orientation,
and let $\overline{\pi}: X\lra Y$ be the collapsing map. If $g$ is
any smooth metric on $Y$ such that $\overline{\pi}^* g$ is
compatible with the complex structure on $X$, then there is a
one-to-one correspondence between
\begin{enumerate}
\item equivalence classes of $g$-self-dual Yang-Mills connections on
a unitary bundle $E_{top}$ over Y, and
\item equivalence classes of holomorphic bundles $E$ on $X$
topologically isomorphic to $\overline{\pi}^* E$ whose restriction
to $L_{\infty}$ is holomorphically trivial and is equipped with a
compatible unitary structure.
\end{enumerate}

\end{theorem}

For the case of $Y=S^4$, i.e. when $n=0$, then $X = \mathbb{C}P^2$,
and $L_\infty$ can be taken to be a line in $\mathbb{C}P^2$.

In view of this result, we will take the necessary steps to
establish that the moduli space of instantons for $S^4$ with a
Bonneau metric is smooth and moreover diffeomorphic to the one for
$S^4$ with a round metric.

Take the round metric given by
$$g=\frac{dr^2+ r^2 + r^2 [ (\sigma^1)^2 + (\sigma^2)^2 + (\sigma^3)^2]}{(1+r^2)^2}$$
where $r$ is a radial coordinate with $r\in (0, +\infty)$. A compatible almost complex structure
is given by (1,0)-forms spanned by
$$\begin{array}{lcl} \eta^1 & = & dr + i r\sigma^3\\ \eta^2 & = & \sigma^1 + i \sigma^2  \end{array}$$
This almost complex structure extends over to $r=0$ and will be denoted by $J_r$.

Now, let us consider the Bonneau metrics. Here a compatible almost
complex structure on $S^4\backslash\{0, \infty\}$ is the one given
by taking the (1,0)-forms to be spanned by
$$\begin{array}{lcl} \theta^1  & = &  a dx + i c \sigma^3\\
\theta^2 & = & \sigma^1 + i \sigma^2 \end{array}$$
As we will see later this structure extends for $x=-\infty$.  For now let us check that
this is actually integrable on $S^4\backslash\{\infty\}$. We have
$$\begin{array}{lcl}
d\theta^1 & = & i \frac{c'}{a}\, \theta^1 \wedge \sigma^3 + i c\, \theta^2 \wedge \sigma^2 \\
d\theta^2 & = & -i \, \theta^2 \wedge \sigma^3
 \end{array}$$ so both $d\theta^1$ and $d\theta^2$ are in the ideal
 generated by $\{ \theta^1, \theta^2\}$. Call this complex structure
 $J_B$.

We wish to construct a diffeomorphism of $S^4$ such that it
sends one almost complex structure into the other. It suffices to
find a coordinate $R\in (0, +\infty)$ such that
$$f (dR + i R\sigma^3) = a dx + i c \sigma^3$$ for some smooth function
$f$. We have that $R$ satisfies the following
$$\left\{ \begin{array}{lcl} f dR = a dx \\ f R = c \end{array} \right.$$
Then $$\frac{dR}{R} = \frac{a}{c}\, dx \quad \Rightarrow  \quad
\log(R) = \int\frac{a}{c}\, dx $$ We now wish to show that this
extends smoothly at $x=k$ and $x = -\infty$. Calculating the
asymptotic expansion around $x=k$ for $\disfrac{a}{c}$, we have
$$\frac{a}{c} = (k-x)^{-1} + O(k-x).$$
Then
$$\frac{a}{c} \sim (k-x)^{-1} \quad \Rightarrow \quad \log(R) \sim -\log(k-x) \quad \Rightarrow R \sim \frac{1}{k-x}.$$
For $x=-\infty$, we have
$$\frac{a}{c} = - x^{-1} + O(x^{-2})$$ and so near $-\infty$, $R \sim
\disfrac{1}{x}$. In particular the complex structure compatible with
the Bonneau metrics extends to $S^4\backslash\{ \infty\}$. We have
then a diffeomorphism
$$\varphi: (S^4\backslash\{ \infty\}, J_r) \lra (S^4\backslash\{ \infty\}, J_B).$$
We consider the twistor space $Z$ to $S^4$ with a Bonneau metric.
The complex structure defined above on $S^4\backslash\{\infty \}$ is
compatible with the metric and gives a section of $$Z\lra
S^4\backslash\{\infty \}.$$ On the other hand, the diffeomorphism
$\varphi$ identifies this with the complex structure of
$\mathbb{C}P^2\backslash\mathbb{C}P^1$. The fact that the
diffeomorphism  $\varphi: S^4 \lra S^4$ commutes with the
$U(2)$-action means that
$$D\varphi_{\infty}:T_{\infty}\lra T_{\infty}$$ is conformal (given by multiplication
by a scalar) and so the $\mathbb{C}P^1$ over $\infty$ is sent to
$L_\infty\subset X$ and the almost complex structures correspond. We
can construct the diagram
$$\begin{CD} \mathbb{C}P^2 @>>> X\\ @V\pi VV @VV\overline{\pi}V \\ S^4 @>\varphi >> S^4 \end{CD}$$
where $X$ is  biholomorphic to $\mathbb{C}P^2$.   Thus
$\overline{\pi}^*(g)$, where $g$ is a Bonneau metric, is compatible
with the complex structure on $X$.

\medskip

\begin{remark}
Note that the mapping $r\lms \disfrac{1}{r}$ provides the same
result for a complex structure with the opposite orientation.
\end{remark}

We have therefore checked all the conditions of theorem \ref{Theo:
buch} and hence we deduce:

\medskip
\begin{theorem}\label{Theo: Bonneau-moduli}
Let $\mathcal{M}_B$ be the moduli space of $SU(2)$-self-dual
connections of charge 1 for a Bonneau metric on $S^4$. Then
$\mathcal{M}_B$ is diffeomorphic to $\mathcal{M}$, the moduli space
of $SU(2)$-self-dual connections of charge 1 for a round metric on
$S^4$.
\end{theorem}

\medskip

\medskip

If we consider the characterization of this moduli space in terms of
$\mathbb{H}^5=\mathbb{R}^4 \times \mathbb{R}^+ = \mathbb{C}^2\times \mathbb{R}^+$
and identify $\mathcal{M}_B$ with $\mathcal{M}$, then the
$U(2)$-invariant instanton equivalence classes are given by the
curve $$\{(z_1, z_2, t): z_1 = z_2 = 0 \}$$
and this contains the equivalence classes of the connections $\nabla^+$ and $\nabla^-$, the connections with
skew torsion $H$ and $-H$ respectively.

An interesting question is whether or not $\nabla^+$ and $\nabla^-$ define the same point on this line, i.e. are gauge equivalent. The answer is no. We need only a counter-example so for simplicity we can choose the parameter $k$ in $(\ref{Eq: Bonneau-metrics})-(\ref{Eq: Bonneau-forms})$ to be zero. We can then proceed to a somewhat lengthy calculation which goes as follows: suppose that there is a $SU(2)$-gauge transformation $g: \S^+ \lra S^+$ such that $g^{-1} \nabla^+ g = \nabla^-$. Then $g$ is a section of $\S^+\otimes S^+$ which is convariantly constant under the tensor product connection $\nabla = \nabla^+\otimes 1 + 1 \otimes \nabla^-$. In this case, $g$ will be annihilated by the curvature of $\nabla$, $R^\nabla$. We can check that there is only one $g$ such that $R^\nabla g = 0$ and that it is totally determined by the metric $ds^2$ and the three form $H$ extending both at $x=k$ and $x=-\infty$. We can then compute $\nabla g$ and see that $g$ is not parallel.

\medskip

\begin{remark}
Using the reduction procedure of Cavalcanti, \cite{Cav}, the closed three form $H$ will induced a closed three form on $\mathcal{M}_B$, giving it the structure of a manifold with skew torsion.
\end{remark}

\bigskip

{\bf Acknowledgments: } The author would like to thank Nigel Hitchin for pointing her to this topic, and for many helpful conversations that have ensued. This research was partially supported by the Portuguese Foundation for Science and Technology (FCT) through the POPH-QREN scholarship program and by the Research Center of Mathematics of the University of Minho through the FCT pluriannual funding program.

\end{document}